\documentstyle[12pt,newlfont]{amsart}

\newcommand{\nc}{\newcommand}

\nc{\one}{\mbox{\bf 1}}
\nc{\invtensor}{\underset{\leftarrow}{\otimes}}
\nc{\ad}{\operatorname{ad}}
\nc{\rk}{\operatorname{rank}}
\nc{\corank}{\operatorname{corank}}
\nc{\Sym}{\operatorname{Sym}}
\nc{\sym}{\operatorname{sym}}
\nc{\id}{\operatorname{id}}
\nc{\htt}{\operatorname{ht}}
\nc{\Ker}{\operatorname{Ker}}
\nc{\im}{\operatorname{Im}}
\nc{\osp}{\operatorname{osp}}
\nc{\ssp}{\operatorname{sp}}
\nc{\sgn}{\operatorname{sgn}}
\nc{\F}{\operatorname{F}}
\nc{\intl}{\operatorname{int}}

\nc{\Inj}{\operatorname{Inj}}
\nc{\Hom}{\operatorname{Hom}}
\nc{\End}{\operatorname{End}}
\nc{\supp}{\operatorname{supp}}
\nc{\Card}{\operatorname{Card}}

\nc{\Ann}{\operatorname{Ann}}
\nc{\Ind}{\operatorname{Ind}}
\nc{\Coind}{\operatorname{Coind}}

\nc{\wt}{\operatorname{wt}}
\nc{\ch}{\operatorname{ch}}
\nc{\Stab}{\operatorname{Stab}}
\nc{\Sch}{{\cal S}\mbox{\em ch}}

\nc{\Spec}{\operatorname{Spec}}
\nc{\Prim}{\operatorname{Prim}}
\nc{\Aut}{\operatorname{Aut}}
\nc{\Fract}{\operatorname{Fract}}
\nc{\gr}{\operatorname{gr}}
\nc{\wdM}{\widetilde{M}}
\nc{\wdV}{\widetilde{V}}
\nc{\wdeps}{\widetilde{\epsilon}}
\nc{\wdC}{\widetilde{\Bbb C}}

\nc{\Dglie}{\operatorname{{\cal D}glie}}
\nc{\Fin}{\operatorname{{\cal F}in}}

\nc{\Sg}{{\cal S}({\frak g})}

\nc{\Ug}{{\cal U}({\frak g})}
\nc{\Zg}{{\cal Z}({\frak g})}
\nc{\tZg}{{\widetilde{\cal Z}({\frak g})}}
\nc{\Zk}{{\cal Z}({\frak k})}
\nc{\Sh}{{\cal S}({\frak h})}
\nc{\Uh}{{\cal U}({\frak h})}

\nc{\Uk}{{\cal U}({\frak k})}
\nc{\Ag}{{\cal A}({\frak g})}

\nc{\cZ}{\cal Z}
\nc{\cS}{\cal S}
\nc{\cP}{\cal P}
\nc{\cL}{\cal L}
\nc{\cU}{\cal U}
\nc{\cH}{\cal H}
\nc{\cK}{\cal K}
\nc{\cF}{\cal F}
\nc{\fg}{\frak g}
\nc{\CO}{\cal O}
\nc{\fn}{\frak n}
\nc{\fm}{\frak m}
\nc{\fh}{\frak h}
\nc{\ft}{\frak t}
\nc{\fk}{\frak k}
\nc{\fp}{\frak p}
\nc{\fI}{\frak I}

\nc{\dirlim}{\underset{\rightarrow}{\lim}\,} 
\nc{\nen}{\newenvironment}
\nc{\ol}{\overline}
\nc{\ul}{\underline}
\nc{\ra}{\rightarrow}
\nc{\lra}{\longrightarrow}
\nc{\Lra}{\Longrightarrow}
\nc{\Lla}{\Longleftarrow}
\nc{\Llra}{\Longleftrightarrow}
\nc{\thla}{\twoheadleftarrow}
\nc{\hra}{\hookrightarrow}
\nc{\iso}{\overset{\sim}{\lra}}
\nc{\ssubset}{\underset{\not=}{\subset}}

\nc{\Thm}[1]{Theorem~\ref{#1}}
\nc{\Prop}[1]{Proposition~\ref{#1}}
\nc{\Lem}[1]{Lemma~\ref{#1}}
\nc{\Cor}[1]{Corollary~\ref{#1}}
\nc{\Conj}[1]{Conjecture~\ref{#1}}
\nc{\Claim}[1]{Claim~\ref{#1}}
\nc{\Defn}[1]{Definition~\ref{#1}}
\nc{\Exa}[1]{Example~\ref{#1}}
\nc{\Rem}[1]{Remark~\ref{#1}}
\nc{\Note}[1]{Note~\ref{#1}}
\nc{\Quest}[1]{Question~\ref{#1}}
\nc{\Hyp}[1]{Hypoth\`ese~\ref{#1}}
\nen{thm}[1]{\label{#1}{\bf Theorem.\ } \em}{}
\nen{prop}[1]{\label{#1}{\bf Proposition.\ } \em}{}
\nen{lem}[1]{\label{#1}{\bf Lemma.\ } \em}{}
\nen{cor}[1]{\label{#1}{\bf Corollary.\ } \em}{}
\nen{conj}[1]{\label{#1}{\bf Conjecture.\ } \em}{}
\nen{claim}[1]{\label{#1}{\bf Claim.\ } \em}{}

\nen{defn}[1]{\label{#1}{\bf Definition.\ } }{}
\nen{exa}[1]{\label{#1}{\bf Example.\ } }{}

\nen{rem}[1]{\label{#1}{\em Remark.\ } }{}
\nen{note}[1]{\label{#1}{\em Note.\ } }{}
\nen{exer}[1]{\label{#1}{\em Exercise.\ } }{}
\nen{sket}[1]{\label{#1}{\em Sketch of proof.\ } }{}
\nen{quest}[1]{\label{#1}{\bf Question.\ } \em}{}
\nen{hyp}[1]{\label{#1}{\bf Hypoth\`ese.\ } \em}{}
\begin{document}

%top matter
\title[]{On the ghost centre of Lie superalgebras }

\author[]{Maria Gorelik}
\address{
{\tt email: gorelik@@msri.org} 
}

\begin{abstract}
We define a notion of {\em ghost centre} of a Lie superalgebra 
$\fg=\fg_0\oplus\fg_1$ which 
is a sum of invariants with respect to the usual adjoint action 
(centre) and invariants with respect to a twisted adjoint action
(``anticentre''). We calculate the anticentre in the case when
the top external degree of $\fg_1$ is a trivial $\fg_0$-module.
We describe the Harish-Chandra 
projection of the ghost centre for basic classical Lie superalgebras
and show that for these cases the ghost centre
coincides with the centralizer of the even part of the enveloping 
algebra.

The ghost centre of a Lie superalgebra plays a role of the
usual centre of a Lie algebra in some problems of representation 
theory. For instance, for $\fg=\osp(1,2l)$ the annihilator of a Verma
module is generated by the intersection with the ghost centre.
\end{abstract}

\thanks{
The author was partially supported by 
 Chateubriand fellowship and TMR Grant No. FMRX-CT97-0100}

\maketitle

\section{Introduction}
\subsection{}
Let $\fg_0$ be a complex finite dimensional  Lie algebra and $\cZ(\fg_0)$
be the centre of its universal enveloping algebra. Then $\cZ(\fg_0)$
acts on a simple $\fg$-module by an infinitesimal character 
and consequently, such characters separate representations.
Moreover, in the case when $\fg_0$ is semisimple, the annihilator
of a Verma module is generated by the kernel of the corresponding
infinitesimal character.

Let $\fg=\fg_0\oplus\fg_1$ be a complex finite dimensional  
Lie superalgebra
and $\Zg$ be the (super)centre of its universal enveloping algebra $\Ug$.
All $\fg$-modules considered below are assumed to be ${\Bbb Z}_2$-graded
and ``$\fg$-simple module'' means simple as graded module.
The centre $\cZ(\fg)$ 
acts on a simple $\fg$-module by an infinitesimal character, but,
even in the ``nice'' case $\fg=\osp(1,2l)$, the annihilator
of a Verma module is not always generated by the kernel of the 
corresponding
infinitesimal character. In~\cite{gl} we described, for the
case $\fg=\osp(1,2l)$, a polynomial subalgebra $\tZg$ of $\Ug$
which acts on a simple module by ``supercharacter''. The annihilator
of a Verma module is generated by the kernel of the corresponding
supercharacter.

In this paper we introduce a notion of {\em ghost centre} 
$\tZg$ (see~\ref{deftzg}).
This is a subalgebra of $\Ug$ which contains both $\Zg$ and
a centre of $\Ug$ considered as  associative algebra.
The algebra $\tZg$ acts on a simple module by ``supercharacter''
and separates graded representations (see~\ref{evencase}). 

By definition, $\tZg$ is a sum of $\Zg$ and 
so-called {\em anticentre} $\Ag$.
The last one is the set of invariants of $\Ug$ with respect to 
a `nonstandard adjoint action'  $\ad'\fg$ introduced
in~\cite{abf}. The product of two elements from the anticentre lies
in the centre and the product of an element from the centre and
an element from the anticentre belongs to the anticentre.

As well as $\Zg$ itself, $\tZg$ is not easy to
describe and, in general, it is not noetherian
algebra. However, in the case when the top external degree 
$\Lambda^{top}\fg_1$ of $\fg_1$ is a trivial $\fg_0$-module,
$\Ag$ itself as well as its image in the symmetric algebra
can be easily described--- see~\Thm{mapphi}. 
The above condition on $\Lambda^{top}\fg_1$ holds for  the
simple finite-dimensional Lie superalgebras apart from the $W(n)$ type.

The existence of non-zero anticentral elements implies two ``negative''
results. The first one is that the direct generalization of the 
Gelfand-Kirillov conjecture does not hold for the Lie superalgebras 
with even dimensional $\fg_1$---see~\ref{gk}. 
The second one is that Separation theorem does not hold 
for the classical basic Lie superalgebras apart from the 
simple Lie algebras and the superalgebras $\osp(1,2l)$---see~\ref{spr}.

\subsection{}
In the case $\fg=\osp(1,2l)$ Arnaudon, Bauer, Frappat (\cite{abf}) 
and Musson (\cite{mu}) constructed a remarkable even element $T$ 
in the enveloping algebra $\Ug$.
This element is $\ad'\fg$-invariant and
its Harish-Chandra projection is the product of hyperplanes
corresponding to the positive odd roots. The element $T$  
has been called `Casimir's ghost' in~\cite{abf}, since its
square belongs to the centre. 

In~\ref{descrAg} we construct such element $T\in\Ag$ for any 
$\fg$ such that $\Lambda^{top}\fg_1$ is a trivial $\fg_0$-module.
The image  of $T$
in the symmetric algebra belongs to $\Lambda^{top}\fg_1$.
In Section~\ref{HCAg} we show that in the case when 
$\fg$ is a basic classical 
Lie superalgebra, the Harish-Chandra projection of
$T$ is also the product of hyperplanes
corresponding to the positive odd roots. 

In~\cite{s2} A.~Sergeev described 
the set of `anti-invariant polynomials' which are the invariants
of the dual algebra $\Ug^*$ with respect to the  nonstandard
adjoint action  $\ad'\fg$.

\subsection{Content of the paper}
In Section~\ref{Vc} we define our main objects:
the anticentre $\Ag$ and the ghost centre $\tZg$. 
We describe the action of $\tZg$
on the modules of finite length in the case when $\fg$ is 
finite dimensional.

In Section~\ref{existT} we show that $\Ag$ is equal to zero if
$\dim\fg_1$ is infinite. Moreover all elements of $\Ag$ are either
even (if $\dim\fg_1$ is even) or odd (otherwise).
We describe the structure of $\Ag$ 
when $\Lambda^{top}\fg_1$ is a trivial $\fg_0$-module. Namely, we prove
that as a vector space $\Ag$ is isomorphic to $\cZ(\fg_0)$.
The central step of the proof is~\Thm{indcoind} which states
that for any $\fg_0$-module $L$, the induced $\fg$-module
$\Ind_{\fg_0}^{\fg} L$ and the coinduced $\fg$-module 
$\Coind_{\fg_0}^{\fg} L$
are isomorphic up to grading. This theorem allows us to define $T$
as a unique up to a scalar $\ad'\fg$-invariant element
inside the $\ad'\fg$-module generated by $1$. 

In Section~\ref{HCAg} we consider
the case when $\fg$ is a complex classical basic Lie superalgebra. 
In this case, the Harish-Chandra projection of 
$\Zg$ is described by Kac and Sergeev (see~\cite{s1}).
In~\Cor{P(A)}, we describe the Harish-Chandra projection of $\Ag$.

We say that an element $u\in\Ug$ acts on a 
module $M$ by a {\em superconstant} 
if it acts by the multiplication by a
scalar on each graded component $M_{\ol i}$ ($i=0,1$). 
In the case when $\fg$ is finite dimensional and $\dim\fg_1$ is even,
any element of $\tZg$ acts  on a 
simple module $M$ by a superconstant (see~\ref{evencase}).
In~\Cor{supercaction} we show that if $\fg$ is a basic classical 
Lie superalgebra then any element of $\Ug$ acting
by a superconstant on each simple finite dimensional module
belongs to $\tZg$. Moreover $\tZg$ coincides with the centre 
(and the centralizer) of the even part $\Ug_{\ol 0}$ of 
the universal enveloping algebra.
For the case $\fg=\osp(1,2l)$ the last result was proven in~\cite{gl}.

\subsection{}
{\em Acknowledgement.}
I wish to express my gratitude to V.~Hinich and A.~Vaintrob for reading
this paper and providing numerous useful suggestions. I would like
to thank  V.~Serganova who pointed out at an error in an earlier version.
I would like to thank also M.~Duflo, A.~Joseph and E.~Lanzmann 
for helpful discussions.
It is finally a great pleasure to thank my hosts at Strasbourg, 
especially P.~Littelmann, whose hospitality and support are greatly 
appreciated.

\section{Ghost centre}
\label{Vc}
In this paper the ground field is ${\Bbb C}$. 
Let $\fg=\fg_0\oplus\fg_1$ be a finite dimensional Lie superalgebra
such that $\fg_1\not=0$. Everywhere except~\ref{maindef},
$\fg$ is assumed to be finite dimensional.
All $\fg$-modules are assumed to be ${\Bbb Z}_2$-graded. We denote
by $\Pi$ the parity change functor: $\Pi(M):=M\otimes \Pi({\Bbb C})$
where $\Pi({\Bbb C})$ is the trivial odd representation.
Denote by $\Ug$ the enveloping superalgebra of $\fg$
and by $\Zg$  the (super)centre of $\Ug$. 

\subsection{}
\label{maindef}
For a homogeneous $u\in\Ug$ denote by $d(u)$ its ${\Bbb Z}_2$-degree.
For a $\Ug$-bimodule $M$ one defines the adjoint action 
of $\fg$ on $M$ by setting 
$(\ad g) m:=gm-(-1)^{d(g)d(m)} mg$ where $m\in M,g\in \fg$ are homogeneous
elements and $d(m)$ denotes the ${\Bbb Z}_2$-degree of $m$.
Define a twisted adjoint action $\ad'$ of $\fg$ on $M$ as the 
adjoint action 
of $\fg$ on the bimodule $\Pi(M)$. One has
$$(\ad'g)(u)=gm-(-1)^{d(g)(d(m)+1)}mg.$$
Assume that  $M$ has a  superalgebra structure such that $g(m_1m_2)
=g(m_1)m_2$ and $(m_1m_2)g=m_1(m_2)g$ for all $g\in\fg,\ m_1,m_2\in M$.
Then for any homogeneous $m_1,m_2\in M$ and $g\in\fg$ one has
$$\begin{array}{ll}
(\ad' g) (m_1m_2)& =((\ad g)m_1)m_2+(-1)^{d(g)d(m_1)}m_1((\ad'g)m_2)\\
 & =((\ad' g)m_1)m_2+(-1)^{d(g)(d(m_1)+1)}m_1((\ad g)m_2).
\end{array}$$
Moreover if  $m$ is $\ad'\fg$-invariant then
\begin{equation}
\label{theta}
\begin{array}{lc}
(\ad' g) (m_1m)=((\ad g)m_1)m, & (\ad g) (m_1m)=((\ad' g)m_1)m.
\end{array}
\end{equation}
\subsubsection{}
\begin{exa}{ex1}
Let $N$ be a $\Ug$-module and $\End(N)$ be the ring 
of its ${\Bbb C}$-linear endomorphisms. Then $\End(N)$ 
admits a natural structure of graded $\Ug$-bimodule. Let $\theta$
be the endomorphism of $N$ which is equal to $\id$ (resp., $-\id$) 
on the even
(resp., odd) component of $N$. Then $\theta$ is an even 
$\ad'\fg$-invariant
homomorphism which commutes with the even elements of $\End(N)$ and
anticommutes with the odd elements of $\End(N)$. 
The formulas~(\ref{theta})
imply that $\End(N)$ considered as $\ad\fg$-module
is isomorphic to $\End(N)$ considered as $\ad'\fg$-module.
The similar assertion fails
for $\Ug$ (the structure of $\Ug$ as $\ad'\fg$-module is given 
in~\Lem{ad1}).
\end{exa}

\subsubsection{}
\label{deftzg}
Let us call {\em anticentre\/} $\Ag$ 
the set of elements of $\Ug$ which are invariant with respect to $\ad'$.
Remark that any even element of the anticentre  anticommutes with
odd elements of  $\Ug$ and commutes with even ones and any odd element
of the anticentre commutes with all elements of  $\Ug$.
Clearly the anticentre is a module over the centre and the product
of any two elements of the anticentre belongs to the centre.
For example, for $\fg=\osp(1,2l)$ $\Ag$ is a free rank one module
over $\Zg$ (see~\cite{gl}, 4.4.1).
This is not true for a general Lie superalgebra.

Let us call {\em  ghost centre\/} $\tZg$ the sum of $\Ag$ and $\Zg$.
It is clear that $\tZg$ is a subalgebra of $\Ug$ which contains
the centre of $\Ug$ considered as associative algebra.

In order to describe the action of $\tZg$ on simple modules,
note the following version of Schur's lemma for Lie superalgebras
\subsubsection{}
\begin{lem}{lemSch}
Let $\fg$ be a finite or countable dimensional Lie superalgebra and
$M=M_{\ol 0}\oplus M_{\ol 1}$ be a simple $\fg$-module.
Then either $\End(M)^{\ad\fg}=k\id$ or
$\End(M)^{\ad\fg}=k\id\oplus k\sigma$ where the odd element $\sigma$ 
provides a $\fg$-isomorphism $M\iso\Pi(M)$ and $\sigma^2=\id$.
\end{lem}
\begin{pf}
Assume that $\phi\in\End(M)^{\ad\fg}$ is even. Both homogeneous components
$M_{\ol 0}$ and $M_{\ol 1}$ are simple $\Ug_{\ol 0}$-modules.
Since $\Ug_{\ol 0}$ is a complex countable dimensional 
associative algebra,
the restriction of $\phi$ on $M_{\ol 0}$ (resp., on $M_{\ol 1}$)
is some constant $c_0$ (resp., $c_1$)--- see~\cite{bz}.
Since $M$ is simple, $c_0=c_1$ and so $\phi=c_0\id$.

Assume that $\phi\in\End(M)^{\ad\fg}$ is odd. Then $\phi^2$ is even
and so $\phi^2=c\id$ for some $c\in {\Bbb C}$.
If $c=0$ then $\Ker\phi\not=0$
and so $\phi=0$. Otherwise $\phi$ is invertible and provides 
a $\fg$-isomorphism $M\iso\Pi(M)$. Set $\sigma=\phi/\sqrt {c}$.
Let $\psi$ be another odd $\ad\fg$-invariant endomorphism
such that $\psi^2=\id$. Then $(\psi\pm\sigma)$ are also 
odd $\ad\fg$-invariant endomorphisms. Therefore $(\psi+\sigma)$ (resp.,
$(\psi-\sigma)$) is either isomorphism or zero. Since
$(\psi+\sigma)(\psi-\sigma)=0$, it implies that $\psi=\pm\sigma$.
This proves the lemma.
\end{pf}

\subsubsection{}
Using~\Exa{ex1}, we conclude that 
$\End(M)^{\ad'\fg}=\End(M)^{\ad\fg}\theta$.
This implies the following lemma describing
the action of $\tZg$ on simple modules.

\begin{lem}{lemacttZg} 
Let $\fg$ be finite or countable dimensional Lie superalgebra,
$M=M_{\ol 0}\oplus M_{\ol 1}$ be a simple $\fg$-module and $z$ 
be an element of $\tZg$. Then the action of $z$ on $M$ is proportional to
$$\begin{array}{ll}
\id, & \text { if } z\in\Zg \text { and } z \text { is even, }\\
0, & \text { if } z\in\Zg \text { and } z \text { is odd, }\\
\theta, & \text { if } z\in\Ag \text { and } z \text { is even, }\\
\sigma\theta, & \text { if } z\in\Ag \text { and } z \text { is odd.}
\end{array}$$
\end{lem}

\subsection{Case $\dim\fg_1$ is even}
\label{evencase}
In this case all elements of $\Ag$ are even (see~\ref{infini}).
Denote by $\wdC$ the algebra spanned by $\id$ and $\theta$. Then
$\wdC={\Bbb C}[\theta]/(\theta^2-1)$. Denote by $\pi$
the algebra involution of $\wdC$ sending $\theta$ to $-\theta$.

\begin{defn}{defsch}
An algebra homomorphism $\chi:\ \tZg\to \wdC$ is called 
{\em supercharacter} if $\chi(\Zg)={\Bbb C}$ and 
$\chi(\Ag)\subseteq {\Bbb C}\theta$.
\end{defn}

By~\Lem{lemacttZg}, $\tZg$ acts on a simple modules $M$ by a 
supercharacter $\chi_M$. Moreover $\chi_{\Pi(M)}=\pi\chi_M$. 

\subsubsection{}
\label{dec}
The standard consequence of Schur's lemma is the following statement.
Any finite length module $M$ has a unique decomposition
into a direct sum of submodules $M_i$ such that, for any
fixed $i$,  all simple subquotients of $M_i$ have the
same infinitesimal character and these characters
are pairwise distinct for different $i$. 
Similarly, one can deduce from~\Lem{lemacttZg},
that  any finite length module $M$ has a unique decomposition
into a direct sum of submodules $M_j$ such that, for any
fixed $j$,  all simple subquotients of $M_j$ have the same supercharacter
and these  supercharacters are pairwise distinct for different $j$.
This new decomposition is a refinement of the previous one.
For example, let $L$ be a simple module such that $\Ag$ does not
lie in $\Ann L$. Then $L$ and $\Pi(L)$ have different supercharacters.
This, for instance, implies that though they have
the same infinitesimal character, there are no non-trivial extensions
of $L$ by $\Pi(L)$.

\subsection{Case $\dim\fg_1$ is odd}

In this case all elements of $\Ag$ are odd (see~\ref{infini}).
Retain notation of~\ref{evencase}.
The algebra spanned by $\id$ and $\sigma\theta$ (see~\Lem{lemSch})
is isomorphic to  $\wdC$.  However if $L$ is a simple module 
such that $aL\not=0$ for some $a\in\Ag$, then the product of $\theta$
and the image of $a$ in $\End L$ provides an isomorphism
$s:L\iso\Pi(L)$. One can choose $a$ such that $s^2=\id$.
There are two possible choices of such $s$ which differ by sign.
As a consequence, in this case, it is more natural to define 
a {\em odd supercharacter} as a pair of homomorphism $(\chi,\pi\chi)$
where $\chi$ satisfies the conditions given in~\Defn{defsch}
and $\pi$ is the involution of $\wdC$ sending $\sigma\theta$ to
$-\sigma\theta$. Observe that if $L\not\cong\Pi(L)$ then $\chi=\pi\chi$.

As in~\ref{dec}, odd supercharacters allows us to construct a 
decomposition of any module of finite length, but, probably, it
always coincides with the decomposition coming from the infinitesimal
characters.

\subsubsection{}
\begin{exa}{ex2}
Let $\fg_1$ be generated by $x$ and $\fg_0$
be generated by $[x,x]$. Then $\Ug={\Bbb C}[x]$, $\Zg={\Bbb C}[x^2]$ and 
$\Ag={\Bbb C}[x^2]x$ is a cyclic $\Zg$-module generated by $x$.
The list of the simple representations of $\fg$ is the following:

a) Two trivial representations (one is even and one is odd). 
The corresponding odd supercharacter sends $\Ag$ to zero.

b) Two-dimensional representations $L(\lambda)$ 
($\lambda\in {\Bbb C}\setminus\{ 0\}$)
spanned by $v$ and $xv$ where $x^2v=\lambda v$. 
The corresponding odd supercharacter sends $x$ to 
$\pm\sqrt{\lambda}\sigma\theta$.
The representations
$L(\lambda)$ and $\Pi(L(\lambda))$ are isomorphic.
\end{exa}

\section{Anticentre $\Ag$}
\label{existT}
Retain notation of Section~\ref{Vc}. Assume that
$\fg=\fg_0\oplus\fg_1$ 
is a Lie superalgebra such that $\fg_1$ is
finite dimensional and $\Lambda^{top}\fg_1$ is a trivial $\fg_0$-module.
In this section we construct 
a linear injective map from the centre $\cZ(\fg_0)$ to the anticentre
$\Ag$---see~\Thm{mapphi}. 
This allows us to describe the image of $\Ag$ in the symmetric
algebra $\Sg$.

\subsection{}
Denote by $\cF$ the canonical filtration of $\Ug$
given by $\cF^k:=\fg^k$. 
Recall that this is an $\ad\fg$-invariant filtration 
and that the associated graded algebra
$\gr_{\cF}\Ug=\Sg$ is supercommutative. For $u\in\Ug$ denote its image in
$\Sg$ by $\gr u$. Remark that $(\ad'x)(u)=2xu-(\ad x)(u)$
for $x\in\fg_1$ and $u\in\Ug$. Therefore 
\begin{equation}
\label{grad}
\gr ((\ad'x)(u))=2(\gr x)(\gr u),\ \ \forall u\in\Ug,x\in\fg_1
\text { s.t. } \gr (xu)=(\gr x)(\gr u).
\end{equation}

\subsubsection{}
Let $L$ be an even vector space endowed by a structure of 
$\fg_0$-module.
Denote by $\Ind_{\fg_0}^{\fg}L$ the supervector space
$\Ug\otimes_{\cU(\fg_0)} L$ (here $\Ug$ is considered as a right
$\cU(\fg_0)$-module) equipped with the natural  left 
$\Ug$-module structure. 

Let $L$ be a submodule of  $\cU(\fg_0)$ with respect to
$\ad\fg_0$-action. Denote by $(\ad' \fg) (L)$ the $\ad'\fg$-submodule
of $\Ug$ generated by $L$. Note that there is a natural surjective map 
from $\Ind_{\fg_0}^{\fg}L$ to 
$(\ad' \fg) (L)$ given by $u\otimes m\mapsto (\ad'u)m$ for 
$u\in\Ug,m\in L$.

\subsubsection{}
\begin{lem}{ad1}
Let $L$ be a submodule of  $\cU(\fg_0)$ with respect to
$\ad\fg_0$-action. The natural map 
$\Ind_{\fg_0}^{\fg}L\to (\ad' \fg) (L)$ is an isomorphism.  Moreover
$\Ug=(\ad'\Ug)\cU(\fg_0)$ and thus as $\ad'\fg$-module
 $\Ug$ is isomorphic to $\Ind_{\fg_0}^{\fg}\cU(\fg_0)$.
\end{lem}
\begin{pf}
Let  $\{x_i\}_{i\in I}$ be an ordered basis of $\fg_1$.
For any finite subset $J\subseteq I$ set  $x_J:=\prod_{i\in J} x_i$, 
where the product is taken with respect
to the fixed order. Then  the elements $\{\gr x_J\}_{J\subseteq I}$
form a basis of $\Lambda \fg_1\subset {\Sg}$. 
Choose a basis $\{u_j\}_{j\in S}$ in $L$ such that 
$\{\gr u_j\}_{j\in S}$ are linearly independent in $\gr_{\cF}\cU(\fg_0)$.
Using~(\ref{grad}) one concludes that 
$\gr (\ad' x_J)u_j=2^{|J|}(\gr x_J)(\gr  u_j)$ for all finite subsets
$J\subseteq I,j\in S$. Therefore the elements
$\{(\ad'x_J)u_j\}_{J\subseteq I,j\in S}$ are linearly independent.
This proves the first assertion. 

For the second assertion,  note that $\gr\Ug$ is spanned by 
the elements of the form $(\gr x_J)(\gr u)$ with $u\in\cU(\fg_0)$.
Now $(\gr x_J)(\gr u)=\gr ((\ad'x_J)u)/2^{|J|}$ and so 
$\gr (\ad'\Ug)\cU(\fg_0)=\gr\Ug$. Therefore $\Ug=(\ad'\Ug)\cU(\fg_0)$
as required.
\end{pf}
The isomorphism $\Ug\cong \Ind_{\fg_0}^{\fg}\cU(\fg_0)$ 
is proven in~\cite{s2}, 3.2. 

\subsubsection{}
\begin{cor}{infini}
If ${\frak g}_1$ has infinite dimension then $\Ag=0$.
If $\dim {\frak g}_1$ is even, all elements of $\Ag$ are even
and if $\dim {\frak g}_1$ is odd, all elements of $\Ag$ are odd.
\end{cor}
\begin{pf}
Retain notation of~\Lem{ad1}.
Any element of  $\cU({\frak g})$ can be written in 
a form $u=\sum_J (ad' x_J) u_J$ 
where $u_J\in \cU(\fg_0)$. Take $u\not=0$ and set 
$m=\max \{|J|\,|\ u_J\not=0\}$. Assume that $m<\dim {\frak g}_1$. Take
$J$ such that  $|J|=m$ and $u_J\not=0$; take $i\in I\setminus J$.
Modulo
$\sum_{|J'|<m+1}  (\ad' x_{J'})\cU(\fg_0)$ one has
$$(ad'x_i)u=\sum_{|J'|=m}(ad' x_i x_{J'}) u_{J'}\not =0.$$ 
Thus if $u\in \Ag$ then $m=\dim {\frak g}_1$. 
Since $\Ag$ is a ${\Bbb Z}_2$-graded
subspace of $\cU({\frak g})$, the assertion follows.
\end{pf}

\subsection{$\Ind$ and $\Coind$}
\label{ico}
Let $L$ be an even vector space endowed by a structure of
left $\fg_0$-module.
Denote by $\Coind_{\fg_0}^{\fg}L$ the supervector space
$\Hom_{\cU(\fg_0)}(\Ug,L)$ (here $\Ug$ is considered as a left
$\cU(\fg_0)$-module) equipped with the following left 
$\Ug$-module structure:
$(uf)(u'):=f(u'u)$ for any $f\in \Hom_{\cU(\fg_0)}(\Ug,L),\ 
u,u'\in \Ug$. The aim of this subsection is to prove that
$\Ind_{\fg_0}^{\fg}L$ and $\Coind_{\fg_0}^{\fg}L$  (resp., 
$\Pi(\Coind_{\fg_0}^{\fg}L)$\,)  are isomorphic
if $\dim\fg_1$ is even (resp., odd).

\subsubsection{}
Retain notation of~\Lem{ad1}. For $k\in {\Bbb N}$ set
$$\cF_o^k:=\sum_{J\subseteq I,|J|\leq k}\cU(\fg_0)x_J.$$
One has $x_Jx_{J'}=\pm x_{J\cup J'}$ modulo $\cF_o^{|J|+|J'|-1}$.
This implies that $\cF_o^p\cF_o^q\subseteq\cF_o^{p+q}$ and
thus $\cF_o$ is a filtration of $\Ug$.
In particular, $\cF_o^k$ are $\cU(\fg_0)$-bimodules and
the filtration does not depend from the choice of $\{x_i\}_{i\in I}$.

Consider $\Ug$ as a left $\cU(\fg_0)$-module.
Denote by $\iota$ a $\fg_0$-homomorphism from $\Ug$ to $\cU(\fg_0)$
such that $\ker \iota=\cF_o^{|I|-1}$
and $\iota(x_I)=1$. Recall that $\ker \iota$ does not depend
from the choice of basis in $\fg_1$. Note that $u_0x_I=x_Iu_0$
modulo $\cF_o^{|I|-1}$ for
any $u_0\in \cU(\fg_0)$, since $\Lambda^{top}\fg_1$
is a trivial $\fg_0$-module. Thus $\iota: \Ug\to\cU(\fg_0)$
is a homomorphism of $\cU(\fg_0)$-bimodules.

Define a map $(.|.)$ from $\Ug\otimes_{\cU(\fg_0)}\Ug$
to $\cU(\fg_0)$ by setting $(u|u')=\iota(uu')$. For any subsets $J,J'$
of $I$ set $\delta_{J,J'}=1$ if $J=J'$ and $\delta_{J,J'}=0$ otherwise.

\subsubsection{}
\begin{lem}{iota}
For any $J\subseteq I$ there exist $u_J,v_J\in\Ug$ 
such that $(u_J|x_{J'})=(x_{J'}|v_{J})=\delta_{J,J'}$.
\end{lem}
\begin{pf}
We prove the existence of $v_J$ by induction on $r=|I\setminus J|$.
For $r=0$, $J=I$ and $v_I=1$ satisfies the conditions.

Fix $J\subseteq I$. For any $J'\subseteq I$ such that $|J'|\leq |J|$,
one has 
$x_{J'}x_{I\setminus J}=\pm x_{I\setminus J\cup J'}$ modulo 
$\ker \iota=\cF_o^{|I|-1}$.
Thus $(x_{J'}|x_{I\setminus J})=0$ for $J\not=J'$ and
$(x_{J}|x_{I\setminus J})=\pm 1$. 
Set
$$v:=x_{I\setminus J}-\sum_{|J'|>|J|}v_{J'}
(x_{J'}|x_{I\setminus J}).$$
Then $(x_{J'}|v)=0$ for any $J'\subseteq I, J'\not=J$ and 
$(x_{J}|v)=\pm 1$. This proves the assertion.

The existence of $u_J$ can be shown similarly.
\end{pf}

\subsubsection{}
\begin{thm}{indcoind}
Assume that $\fg=\fg_0\oplus\fg_1$ is a Lie superalgebra
such that $\Lambda^{top}\fg_1$ is a trivial $\fg_0$-module.
Then for  any  $\fg_0$-module $L$ the linear map $\Psi$  defined by
$$\Psi(u'\otimes m)(u):=(u|u')m,\ \forall m\in L,\ u,u'\in \Ug$$
provides an isomorphism 
$\Ind_{\fg_0}^{\fg}L\iso\Pi^{\dim\fg_1}(\Coind_{\fg_0}^{\fg}L).$
\end{thm}
\begin{pf}
For any $u_0\in\cU(\fg_0)$ one has
$(u|u'u_0)m=\iota(uu'u_0)m=\iota(uu')u_0m=(u|u')u_0m$ and thus
$\phi$ is well-defined on $\Ind_{\fg_0}^{\fg}L=\Ug\otimes_{\cU(\fg_0)} L$.
Moreover
$$\Psi(u'\otimes m)(u_0u)=(u_0u|u')m=u_0(u|u')m=u_0\Psi(u'\otimes m)(u)$$
and so $\Psi(u'\otimes m)$ is a $\cU(\fg_0)$-linear map.

For any $s\in \cU$ one has 
$$\Psi(su'\otimes m)(u)=(u|su')m=(us|u')m=\Psi(u'\otimes m)(us)=
(s\Psi(u'\otimes m))(u)$$
and so $\Psi$ is a homomorphism of left $\Ug$-modules.

Since $\Psi(1\otimes m)(x_J)=\delta_{I,J}m$, the element 
$\Psi(1\otimes m)$
is even iff $x_I\in\Ug$ is even. Thus $\Psi$ is a map 
from $\Ind_{\fg_0}^{\fg}L$ to $\Coind_{\fg_0}^{\fg}L$ if $\dim\fg_1$ 
is even
and from $\Ind_{\fg_0}^{\fg}L$ to $\Pi(\Coind_{\fg_0}^{\fg}L)$ 
if $\dim\fg_1$ is odd.

Any element of $\Ind_{\fg_0}^{\fg}L$
can be written in the form $\sum_{J\subseteq I} x_J\otimes m_J$
where $m_J\in L$. Fix $J'\subseteq I$ and choose $u_{J'}\in\Ug$ 
as in~\Lem{iota}. Then 
$\Psi(\sum_{J\subseteq I} x_J\otimes m_J)(u_{J'})=m_{J'}$.
This implies that $\ker\Psi=0$.
 
Fix $J\subseteq I$ and choose $v_{J}\in\Ug$ 
as in~\Lem{iota}. Then for any $m\in L$ one has
$\Psi(v_{J}\otimes m)(x_{J'})=\delta_{J,J'}m$. This implies the 
surjectivity of $\Psi$ and completes the proof.
\end{pf}

\subsection{}
\label{descrAg}
Retain notation of~\Lem{iota}.

\begin{thm}{mapphi}
Assume that $\fg=\fg_0\oplus\fg_1$ is a Lie superalgebra
such that $\Lambda^{top}\fg_1$ is a trivial $\fg_0$-module.
Then the map $\phi: z\mapsto (\ad' v_{\emptyset})z$ provides 
a linear isomorphism  $\cZ(\fg_0)\iso \Ag$. 
Moreover one has $\gr\phi(z)=x\gr z$ where 
$x$ is an element of $\Lambda^{top}(\fg_1)$.
\end{thm}
\begin{pf}
The proof follows from~\Lem{ad1} and~\Thm{indcoind}. We give a full
detail below.

Denote by $\wdeps$ a trivial even representation of $\fg$ and let $e$
be a non-zero vector of $\wdeps$. There is a canonical 
bijection $\Phi$ from $\Hom_{\fg_0}(\wdeps,\cU(\fg_0))$ onto 
$\Hom_{\fg}(\wdeps,\Coind_{\fg_0}^{\fg}\cU(\fg_0))$ given by
$$\Phi(\psi)(e)(u)=\psi(ue)\ \ \forall \psi\in 
\Hom_{\fg_0}(\wdeps,\cU(\fg_0)),
u\in\Ug.$$
Combining the map $\Phi$ with the natural bijection 
$\cZ(\fg_0)\iso\Hom_{\fg_0}(\wdeps,\cU(\fg_0))$, we obtain 
the bijection 
$$\Phi': \cZ(\fg_0)\iso\Hom_{\fg}(\Pi^{\dim\fg_1}(\wdeps),
\Coind_{\fg_0}^{\fg}\cU(\fg_0))$$ given by
$$(\Phi'(z)(e))(x_J)=\delta_{J,\emptyset}z.$$
In view of~\Thm{indcoind}, $\Phi'$ induces the bijection 
$$\Phi'': \cZ(\fg_0)\iso
\Hom_{\fg}(\Pi^{\dim\fg_1}(\wdeps),\Ind_{\fg_0}^{\fg}\cU(\fg_0))$$
given by 
$$\Phi''(z)(e)=v_{\emptyset}\otimes z.$$ 
Finally,~\Lem{ad1} implies
that the map sending $z$ to $(\ad' v_{\emptyset})z$ provides
a linear isomorphism  $\cZ(\fg_0)\iso \Ag$ and moreover $\Ag$ lies
in $\Ug_{\ol 0}$ if $\dim\fg_1$ is even and in $\Ug_{\ol 1}$
if $\dim\fg_1$ is odd.

The proof of~\Lem{iota} shows that $v_{\emptyset}=\pm x_{I}
+\sum_{J\not=\emptyset} x_{I\setminus J} d_J$ where $d_J$ 
are some elements of $\cU(\fg_0)$. Therefore 
\begin{equation}\label{r2}
\phi(z):=(\ad'v_{\emptyset})(z)=(\ad'(x_I+\sum_{J\ssubset I}c_Jx_J))z
\end{equation}
where $c_J$ are some scalars. By the formula~(\ref{grad}),
$\gr\phi(z)=x\gr z$ for $x:=\gr x_I\in \Lambda^{top}(\fg_1)$.
This completes the proof.
\end{pf}
\subsection{}
The condition on  $\Lambda^{top}(\fg_1)$ is essential for 
both~\Thm{indcoind} and~\Thm{mapphi}. 
\subsubsection{}
\begin{exa}{exaW}
Let $\fg$ be the Lie superalgebra $W(1)$ spanned by the 
even element $g$ and the odd element $x$ subject to 
the relations $[x,x]=0, [g,x]=x$.
Then $\Zg=\tZg={\Bbb C}$ and $\Ag=0$.
\end{exa}

\subsubsection{}
For~\Thm{indcoind} the condition on $\Lambda^{top}$
is not only sufficient but also necessary. In fact, 
assume that the isomorphism $\Psi$ exists for a trivial $\fg_0$-module
$\epsilon$. Then one has the following isomorphisms of $\fg_0$-modules
$$\Lambda\fg_1\iso\Ind_{\fg_0}^{\fg}\epsilon\iso
\Coind_{\fg_0}^{\fg}\epsilon\iso(\Lambda\fg_1)^*.$$
Thus $\Lambda\fg_1$ is a self-dual $\fg_0$-module.
The last is equivalent
to the condition  $\Lambda^{top}(\fg_1)\iso\epsilon$
if $\fg_1$ is assumed to be finite dimensional.
\subsubsection{}
We describe below for which simple Lie superalgebras the condition 
$\Lambda^{top}(\fg_1)\iso\epsilon$ holds.

A classification theorem of Kac (see~\cite{k2}, 4.2.1)
states that any complex simple finite dimensional Lie superalgebra 
is isomorphic either to one of the classical Lie superalgebra or to one
of the Cartan Lie superalgebras $W(n),S(n), {\widetilde S}(n), H(n)$. 

Evidently the condition holds if $\fg_0$ is semisimple
or if $\fg_1\cong\fg_1^*$ as $\fg_0$-module. In particular the condition
holds for all classical  Lie superalgebras. It is easy to check that
the above condition holds also for the Cartan Lie 
superalgebras $S(n), {\widetilde S}(n), H(n)$ and does not 
hold for the Cartan Lie superalgebras $W(n)$ with $n\not=2$.

\subsubsection{}
Retain the assumption of~\Thm{mapphi}.

\begin{defn}{defT}
Denote by $T$ a non-zero $\ad'\fg$-invariant element belonging to
$(\ad'\Ug) (1)$.
\end{defn}

The element $T$ is defined up to a non-zero scalar and it is
even iff $\dim\fg_1$ is even. Observe that, up to a scalar, $T$
is a unique element of the anticentre whose image in $\Sg$
belongs to $\Lambda^{top}(\fg_1)$.

\subsection{Remarks}
\subsubsection{}
Here we consider $\Ug$ as an associative algebra and denote its centre by
$Z$. Evidently $Z\cap\Ug_{\ol 0}=\Zg\cap\Ug_{\ol 0}$ and
$Z\cap\Ug_{\ol 1}=\Ag\cap\Ug_{\ol 1}$. Hence, in the case when $\fg$
is finite dimensional, 
$$\begin{array}{ll}
Z=\Zg\cap\Ug_{\ol 0} & \text{ if }\ \dim\fg_1\ \text{ is even},\\
Z=(\Zg\cap\Ug_{\ol 0})\oplus\Ag & \text{ if }\ \dim\fg_1\ \text{ is odd}.
\end{array}$$
\subsubsection{}
\label{gk}
In most of the cases $\Ug$ is not a domain--- see~\cite{al}.
However, even if $\Ug$ is a domain (for example $\fg=\osp(1,2l)$)
the direct generalization of the Gelfand-Kirillov conjecture does not
hold for Lie superalgebras.

In fact, let $k$ be a field of characteristic zero and 
$A_n(k)$ be a Weyl algebra over $k$. Recall that the centre
of a Weyl skew field $W_n(k)$ coincides with $k$ and
that $A_n({\ol k})=A_n(k)\otimes_k{\ol k}$. Therefore
a Weyl skew field
does not contain non-central elements whose squares are central.
Take any non-zero $a\in\Ag$. If $\dim\fg_1$ is even then $a\not\in Z$, but
$a^2\in Z$. 
This implies that a Weyl skew field and a skew field 
of fractions of $\Ug$ are not isomorphic if $\dim\fg_1$ is even
and $\Lambda^{top}\fg_1$ is a trivial $\fg$-module.

\section{The case of basic classical Lie superalgebras}
\label{HCAg}
In this section $\fg$ denotes a basic classical Lie superalgebra
(see~\cite{k} and~\ref{ntt} below) such that $\fg_1\not=0$.
In this case the dimension of $\fg_1$ is even and so all elements of $\Ag$
 even. In particular, they anticommute with the odd elements of $\Ug$
and commute with the even ones. 

In this section we  show that 
the restriction of the Harish-Chandra projection $\cP$
on $\Ag$ is an injection and  describe its image. 
We also prove that $\tZg$ coincides with the centralizer 
of $\Ug_{\ol 0}$ and with the set of the elements of $\Ug$ acting 
by superconstants on each simple finite dimensional module.

\subsection{Notation}
\label{ntt}
A finite dimensional simple Lie superalgebra $\fg$ is called 
basic classical if $\fg_0$ is reductive and 
$\fg$ admits a non-degenerate invariant bilinear form.
The list of basic classical Lie superalgebras is the following 
as determined by Kac (see~\cite{k}):
$$\begin{array}{ll}
a) &  \text{simple Lie agebras}\\
b) & A(m,n),\ B(m,n),\ C(n),\ D(m,n),\ D(2,1,\alpha),\ F(4),\ G(3).
\end{array}$$

Fix a Cartan subalgebra
${\frak h}$ in ${\frak g}_0$ and a triangular decomposition
$\fg=\fn^-\oplus\fh\oplus\fn^+$. For a $\Ug$-module $M$ and an element
 $\mu\in\fh^*$ set
$M|_{\mu}=\{m\in M|\ hm=\mu(h)m,\ \forall h\in\fh\}$. 
When we use the notation $\Ug|_{\mu}$, the action of $\fg$ on
$\Ug$ is assumed to be the adjoint action. For $\mu\in\fh^*$ 
we say that  $\mu$ is {\em even\/} if $\Ug|_{\mu}$ is a non-zero 
subspace of the even part of $\Ug$. We say that  $\mu$ is 
{\em odd\/}  if $\Ug|_{\mu}$ is a non-zero subspace of
 the odd part of $\Ug$. Since $\fg$ is a basic classical Lie superalgebra,
$\mu$ is either even or odd in the case when $\Ug|_{\mu}\not=0$.

The Harish-Chandra projection $\cP:\Ug\to\Sh$ is
the projection with respect to the following triangular decomposition
$\Ug=\Uh\oplus(\Ug\fn^++\fn^-\Ug)$ (we identify $\Sh$ and $\Uh$).
An element $a$ of $\Ug|_0$ acts on a primitive vector of weight $\mu$
($\mu\in\fh^*$) by  multiplication by the scalar $\cP(a)(\mu)$. 
Thus the restriction of
$\cP$ on $\Ug|_0=\Ug^{\fh}$ is an algebra homomorphism from
$\Ug|_0$ to $\Sh$.

Denote by $\Delta_0$ the set of non-zero even 
roots of ${\frak g}$. Denote by $\Delta_1$ the set of odd 
roots of ${\frak g}$. Set $\Delta=\Delta_0\cup\Delta_1$.
Set
$$\overline{\Delta}_0\!\!:=
\{\alpha\in\Delta_0|\ \alpha/2\not\in\Delta_1\},\ \ \
\overline{\Delta}_1\!\!:=\{\beta\in\Delta_1|\ 2\beta\not\in\Delta_0\}.$$
Note that $\overline{\Delta}_1$ is the set of isotropic roots.
Denote by $\Delta^{+}$ the set of positive roots and define 
$\Delta^{-}, \Delta_0^{\pm},\Delta_1^{\pm},\overline{\Delta}_0^{\pm}$
as usual.

Denote by $W$ the Weyl group of $\Delta_0$.
For any $\alpha\in \Delta_0$ let  $s_{\alpha}\in W$ be the corresponding
reflection. Let $W'$ be the subgroup of $W$
generated by the reflections $s_{\alpha}, \alpha\in \overline{\Delta}_0$.
Note that $W=W'$ iff all odd roots are isotropic. Otherwise (if
$\fg$ is of the type $B(m,n)$ or $G(3)$) $W'$ is 
a subgroup of index two.

Set 
$$\begin{array}{ccc}
\rho_0\!:={1\over 2}\displaystyle\mathop{\sum}_
{\alpha\in \Delta_0^+}\alpha,& 
\rho_1\!:={1\over 2}\displaystyle\mathop{\sum}_
{\alpha\in \Delta_1^+}\alpha,& 
\rho\!:=\rho_0-\rho_1.
\end{array}$$ 
Define the translated action of $W$ on ${\frak h}^*$ by the formula: 
$$w.\lambda=w(\lambda+\rho)-\rho,\ \ \forall \lambda\in 
{\frak h}^*,w\in W.$$
Define the left translated action of $W$ on $\Sh$ by setting
$w.f(\lambda)=f(w^{-1}.\lambda)$ for any $\lambda\in\fh^*$.

Denote by $(-,-)$ a non-degenerate
$W$-invariant bilinear form on ${\frak h}^*$.
\subsubsection{}
\label{grVerma}
For $\lambda\in\fh^*$ denote a graded $\fg$-Verma module of
the highest weight $\lambda$ by $\wdM(\lambda)$ where the grading 
is fixed in such a way that a highest weight vector has degree zero.
By~\cite{lm}, an element of $\Ug$ annihilating the modules 
$\wdM(\lambda)$ for $\lambda$ running through a Zariski dense
subset of $\fh^*$, is equal to zero.

We use the following result which is
a consequence of a theorem of Kac (see~\cite{ja}, 2.4)

\subsubsection{}
\begin{lem}{irrVerma}
Assume that a pair $(n,\alpha)$ belongs to the following set
$$({\Bbb N}^+\times\overline{\Delta}_0^+)\cup
((1+2{\Bbb N})\times(\Delta_1^+\setminus \overline{\Delta}_1^+))$$
and $\lambda\in\fh^*$ is such that $(\lambda+\rho,\alpha)=n$.
Then $\wdM(\lambda)$ contains a primitive vector of the weight
$\lambda-n\alpha$.
If $\alpha\in\overline{\Delta}_1^+$ and  $(\lambda+\rho,\alpha)=0$
then $\wdM(\lambda)$ contains a primitive vector of the weight
$\lambda-\alpha$.
\end{lem}

\begin{rem}{}
Note that the formula for Shapovalov determinants presented
in~\cite{ja}, 2.4 contains misprints;
the correct formula reads as follows
$$\begin{array}{rl}
\det S_{\eta}& = A\cdot B\cdot C, \\
A&=\displaystyle\prod_{n=1}^{\infty} 
\displaystyle\prod_{\gamma\in {\overline{\Delta}_0^+}} 
(h_{\gamma}+(\rho,\gamma)-n(\gamma,\gamma)/2)^{P(\eta-n\gamma)},\\
B&= \displaystyle\prod_{n=1}^{\infty}
\displaystyle\prod_{\gamma\in \Delta_1^+\setminus \overline{\Delta}_1^+}
(h_{\gamma}+(\rho,\gamma)-
(2n-1)(\gamma,\gamma)/2)^{P(\eta-(2n-1)\gamma)},\\
C&=\displaystyle\prod_{\gamma\in {\overline{\Delta}_1^+}}
(h_{\gamma}+(\rho,\gamma))^{P_{\gamma}(\eta-\gamma)}.
\end{array}$$
This formula  can be proven following the proof of~\cite{j}, 6.1---6.11.
\end{rem}

\subsection{}
For $\beta\in\fh^*$ denote by $h_{\beta}$ the element of $\fh$
such that $\mu(h_{\beta})=(\mu,\beta)$ for any $\mu\in\fh^*$.
Set 
$$t:=\prod_{\alpha\in {\Delta}_1^+}(h_{\alpha}+(\rho,\alpha)).$$

\subsubsection{}
\label{WinvofPA}
\begin{lem}{Pa}
The restriction of the Harish-Chandra projection $\cP$ provides
a linear injective map  $\Ag\to t \Sh^{W.}$.
\end{lem}
\begin{pf}
Recall that $a\in\Ug|_0$ acts on a primitive vector of weight $\mu$
($\mu\in\fh^*$) by  multiplication by the scalar $\cP(a)(\mu)$. 
Fix a non-zero $a\in\Ag$. Since $\Ag\subset\Ug|_0$, 
$a$ acts on the even component of $\wdM(\lambda)$
by multiplication by the scalar $\cP(a)(\lambda)$ and
on the odd component of $\wdM(\lambda)$
by multiplication by the opposite scalar.
The intersection of the annihilators of all Verma modules is zero
(see~\ref{grVerma}) and so 
$\cP(a)$ is a non-zero polynomial in $\Sh$.

Choose a pair $(n,\alpha)$ and an element $\lambda$
satisfying the assumption of~\Lem{irrVerma}. Note that $n\alpha$ is even
iff $\alpha$ is even. This implies that 
$\cP(a)(\lambda)=\cP(a)(\lambda-n\alpha)$ for $\alpha\in 
\overline{\Delta}_0^+$
and $\cP(a)(\lambda)=-\cP(a)(\lambda-n\alpha)$ for 
$\alpha\in \Delta_1^+\setminus \overline{\Delta}_1^+$.
Observe that  for  $\alpha\in \overline{\Delta}_0^+\ $ 
$\lambda-n\alpha=s_{\alpha}.\lambda$ and so 
$\cP(a)(\lambda)=\cP(a)(s_{\alpha}.\lambda)$. For fixed 
$\alpha\in \overline{\Delta}_0^+$ 
the set of $\lambda$ such that $(\lambda+\rho,\alpha)\in {\Bbb N}^+$
is a Zariski dense subset of $\fh^*$. Thus $\cP(a)\in\Sh$ 
is $W'.$-invariant. 

Take
$\alpha\in(\Delta_1^+\setminus \overline{\Delta}_1^+)$. 
Then $2\alpha\in \Delta_0^+$ and $\lambda-n\alpha=s_{2\alpha}.\lambda$;
arguing as above we obtain  that $s_{2\alpha}.\cP(a)=-\cP(a)$.
In particular, $\cP(a)$ is divisible by $(h_{\alpha}+(\rho,\alpha))$.

Now take $\alpha\in\overline{\Delta}_1^+$. Then
$\alpha$ is isotropic. In particular, $(\lambda+\rho,\alpha)=0$
implies $(\lambda-\alpha+\rho,\alpha)=0$. 
Using~\Lem{irrVerma}, we conclude that $(\lambda+\rho,\alpha)=0$ implies that
$\cP(a)(\lambda)=(-1)^n\cP(a)(\lambda-n\alpha)$ for any $n\in{\Bbb N}$. 
Therefore $\cP(a)(\lambda)=0$ if $(\lambda+\rho,\alpha)=0$.
Thus $\cP(a)$ is divisible by $(h_{\alpha}+(\rho,\alpha))$.

Hence $\cP(a)$ is divisible by $(h_{\alpha}+(\rho,\alpha))$
for any $\alpha\in\Delta_1^+$. This implies that $\cP(a)$ is divisible by
$t$. Since $t$ is $W'.$-invariant, $\cP(a)/t$ is also $W'.$-invariant.
For any $\alpha\in(\Delta_0^+\setminus\overline{\Delta}_0^+$,
both $\cP(a)$ and $t$ are antiinvariant with respect to the action of 
$s_{\alpha}$. Thus  $\cP(a)/t$ is invariant with respect to the action of 
$s_{\alpha}$ and so $\cP(a)/t$ is $W.$-invariant. 
This completes the proof.
\end{pf}

\subsubsection{}
\label{degP(a)}
Define a filtration on $\fg$  by setting $\cF_u^{0}=0,\ 
\cF_u^{1}=\fg_1,\ \cF_u^{2}=\fg$ and extend it canonically to an 
increasing
filtration on $\Ug$. Let $z\in\cZ(\fg_0)$ have a degree $r$ with respect
to the canonical filtration. Then,
by~(\ref{r2}), $\phi(z)\in \cF_u^{\dim\fg_1+2r}$ and so $\cP(\phi(z))$
is a polynomial of degree less than or equal to 
$(\dim\fg_1+2r)/2=|\Delta_1^+|+r$. 
In particular, $\cP(\phi(1))$ is a polynomial of degree less than or 
equal to $|\Delta_1^+|$ and so it is equal to $t$ up to a non-zero 
scalar. Recall that the map $\phi$ depends on the choice of basis 
$\{x_i\}_{i\in I}$;
choose a basis $\{x_i\}_{i\in I}$ such that $\cP(\phi(1))=t$.

\subsubsection{}
Fix $r\in{\Bbb N}$ and set $Z_r:=\cZ(\fg_0)\cap\cF^r$.
Denote by $S_r$ the space of $W.$-invariant polynomials
of degree less than or equal to $r$. Take $z\in Z_r$.
Combining~\Lem{Pa} and~\ref{degP(a)}, we conclude
that  $(\cP(\phi(z_r))/t)\in S_r$.
Recall that $\gr\cZ(\fg_0)\iso\Sh^W$ 
as graded algebras and so $\dim Z_r=\dim S_r$. 
Since $\phi$ is a linear isomorphism, it follows that 
$\cP(Z_r)=tS_r$. 

\subsubsection{}
\begin{cor}{P(A)}
The restriction of the Harish-Chandra projection $\cP$ provides
a linear bijective map  $\Ag\to t \Sh^{W.}$. 
In particular, $\cP(T)=t$.
\end{cor}

\subsubsection{}
\begin{lem}{nonzerodiv}
Any non-zero element $z\in\Ag$ is a non-zero divisor in $\Ug$.
\end{lem}
\begin{pf}
Assume that $zu=0$. Recall that $z$ acts by multiplication by
$\cP(z)(\lambda)$ (resp., $-\cP(z)(\lambda)$) on the even (resp., odd)
graded component of $\wdM(\lambda)$. Therefore $u$ annihilates $\wdM(\lambda)$
when $\lambda$ is such that $\cP(z)(\lambda)\not=0$.
Since $\cP(z)\not=0$, the set $\{\lambda|\ \cP(z)(\lambda)\not=0\}$
is a Zariski dense subset of $\fh^*$. 
By~\ref{grVerma}, it implies that $u=0$.
\end{pf}

\subsubsection{}
\begin{rem}{}
On the contrary to the central elements, $\gr z$ is a zero-divisor
for any $z\in\Ag$. In fact, by~(\ref{theta}) $(\ad \Ug)(z)=(\ad'\Ug)(1) z$
and thus  $(\ad\Ug)z$ contains $Tz$. Therefore $z$ and $Tz$ have the same
degree with respect to the canonical filtration and so $\gr T\gr z=0$.
In particular, $T^2$ is a central element whose degree is equal to 
$\dim\fg_1$.
\end{rem}

\subsubsection{}
\begin{cor}{}
$$\Zg\cap\Ag=0.$$
\end{cor}
\begin{pf}
For any $z\in \Zg\cap\Ag$ and any odd element $u\in\fg_1$  one
has $zu=0$. Hence $z=0$ by~\Lem{nonzerodiv}.
\end{pf}

\subsection{The structure of $\tZg$}
The algebra $\tZg$ has the following easy realization.
Consider  the algebra $\widetilde{\Sh}:=\Sh[\xi]/(\xi^2-1)$.
Define a map $\cP': \tZg\to \widetilde{\Sh}$
by setting  $\cP'(z)=\cP(z)$ for $z\in\Zg$ and  $\cP'(z)=\cP(z)\xi$ for
$z\in\Ag$. Since $\tZg\subset\Ug|_0$, the restriction of  
$\cP$ on $\tZg$ is an algebra homomorphism. 
Taking into account~\Cor{P(A)}, 
we conclude that $\cP'$ provides an algebra isomorphism
from $\tZg$ onto the subalgebra $(\cP(\Zg)\oplus t \Sh^{W.}\xi)$ of 
$\widetilde{\Sh}$.
 
\subsubsection{}
Assume that  $\fg$ is of the type $B(m,n)$ or $G(3)$.
Then $W'\not=W$ and so $t$ is not $W.$-invariant. Therefore 
$\cP(\Ag)\cap \cP(\Zg)=\{0\}$. Then, using~\Cor{P(A)},
we conclude that the restriction of the Harish-Chandra projection provides
an algebra isomorphism $\tZg\cong (\cP(\Zg)\oplus t \Sh^{W.})$.

In all other cases, $\cP(\Ag)\subset \cP(\Zg)$.

\subsubsection{}
In the case when $\fg=\osp(1,2l)$, $\Zg$ is a polynomial algebra  and 
$\Ag$ is a cyclic $\Zg$ module
generated by $T$. In other cases (when $\fg$ is basic classical
Lie superalgebra) this does not hold. However, a similar
result hold after a certain localization. 

More precisely, if $\fg\not=\osp(1,2l)$ (that is $\fg$
is not of the type $B(0,l)$) then $\cP(\Zg)$ is strictly contained in
$\Sh^{W.}$. However, since the product of two elements from the anticentre
belongs to the centre, $\cP(\Zg)$ contains $t^2\Sh^{W.}$. Set $Q:=T^2$,
$q:=t^2$. Then the localized algebras $\Zg[Q^{-1}]$  and 
$\Sh^{W.}[q^{-1}]$  are isomorphic.
Moreover $\tZg[Q^{-1}]=\Zg[Q^{-1}]\oplus\Ag[Q^{-1}]$ and $\Ag[Q^{-1}]$
is a cyclic $\Zg[Q^{-1}]$-module generated by $T$.

\subsection{The action of $\tZg$ on the simple modules.}
Let us say that an element $u\in\Ug$ acts on a $\Ug$-module $M$
by a superconstant if it acts by a multiplication by a scalar 
on each graded component of $M$. By~\ref{evencase},  each element of 
$\tZg$ acts by a superconstant on any simple module. 
In this subsection we shall prove that actually
$\tZg$ coincides with the set of elements of $\Ug$ which act
by superconstants on each simple finite dimensional module.
Moreover $\tZg$ coincides with the centralizer
of the even part $\Ug_{\ol 0}$ in $\Ug$.

\subsubsection{}
By definition, $\tZg$ lies in the centralizer of $\Ug_{\ol 0}$ in 
$\Ug$ and even in the centre of $\Ug_{\ol 0}$ since  all elements 
of $\tZg$ are even.

Let $A$ be a centralizer of $\Ug_{\ol 0}$ in $\Ug$ and 
$a$ be an element of $A$.
Clearly, $a$ acts by a superconstant on any Verma module.
On the even component of $\wdM(\lambda)$
$a$  acts by $\cP(a)(\lambda)$. Let $f(a)$ be the function $\fh^*\to k$
such that $a$ acts by $f(a)(\lambda)$ on the odd component of 
$\wdM(\lambda)$.

\subsubsection{}
\begin{lem}{f(a)}
For any $a\in A$ the function $f(a):\fh^*\to k$ is polynomial.
\end{lem}
\begin{pf}
Choose $y\in\fn_1^-$ and
$x\in\fn_1^+$ such that $h:=[y,x]\in\Sh$ and $h\not=0$.
For each $\mu\in\fh^*$ choose 
a highest weight vector $v_{\mu}\in \wdM(\mu)$.
Then $yv_{\mu}$ is  odd and so
$$xay v_{\mu}=f(a)(\mu)xyv_{\mu}=f(a)(\mu)h(\mu) v_{\mu}.$$
Since $xay\in\Ug|_0$ one has 
$xay v_{\mu}=\cP(xay)(\mu)v_{\mu}$.
Thus 
\begin{equation}\label{cmu}
f(a)(\mu)h(\mu)=\cP(xay)(\mu).
\end{equation}
This implies that $\cP(xay)(\mu)$ vanishes on the
whole hyperplane $\{\mu|\ h(\mu)=0\}$. Therefore
$h$ divides $\cP(xay)(\mu)$
and so $f(a)=\cP(xay)/h$ is a polynomial.
\end{pf}

\subsubsection{}
 \Lem{f(a)} implies that an element $a\in A$ acts by $\cP(a)(\lambda)$ 
on the even component of $\wdM(\lambda)$ and by
$f(\lambda)$ on the odd component of $\wdM(\lambda)$.
Arguing as in~\ref{WinvofPA}, 
we obtain that $P':=\cP(a)-f(a)$ belongs to $t\Sh^{W.}$.
Similarly, $P:=\cP(a)+f(a)$ belongs to $\Sh^{W.}$
and moreover for any $\alpha\in{\overline{\Delta_1}}^+$ 
one has $P(\lambda-\alpha)=P(\lambda)$ if $(\lambda+\rho,\alpha)=0$.
By~\cite{s1} and~\Cor{P(A)}, this implies that $P=\cP(z)$ 
for some $z\in\Zg$ and $P'=\cP(z')$ for some $z'\in\Ag$. Then 
$a-(z+z')/2$ kills any Verma module and so $a=(z+z')/2$.
This proves that  $\tZg=A$.

The intersection of the annihilators
of all simple highest weight modules is equal to zero (see~\ref{grVerma}).
This implies that the set of elements of $\Ug$ acting
by superconstants on each graded simple finite dimensional module
coincides with $A$. Hence we obtain
\subsubsection{}
\begin{cor}{supercaction}
If $\fg$ is a basic classical Lie superalgebra then the following algebras
coincide
 
i) The algebra of elements of $\Ug$ which act
by superconstants on each graded simple finite dimensional module.

ii) The algebra $\tZg$.

iii) The centre of $\Ug_{\ol 0}$.

iv)  The centralizer of $\Ug_{\ol 0}$.
\end{cor}

\subsection{A remark concerning the separation theorem.}
\label{spr}
An important structure theorem of Kostant asserts that for
any finite dimensional semisimple Lie algebra 
there exists an $\ad\fg$-submodule $\cH$ of $\Ug$
such that the multiplication map induces the isomorphism
$\cH\otimes \Zg\iso\Ug$. In~\cite{mu}, I.~Musson proved
the similar assertion for $\fg=\osp(1,2l)$. 
These theorems are called the separation theorems. 
We shall show that Separation theorem does not hold for 
any basic classical Lie superalgebra apart from finite dimensional
simple Lie algebras and $\fg=\osp(1,2l)$.

\subsubsection{}
Denote by $V$ (resp., $\wdV$) a trivial representation of $\fg_0$ 
(resp., $\fg$). 

\begin{lem}{lemsepthm}
Assume that $\fg_1$ contains a non-zero element $x$ such
that $[x,x]=0$. Then the trivial submodule $\wdV$ is not a direct summand
of $\Ind_{\fg_0}^{\fg} V$.
\end{lem}
\begin{pf}
Denote a generator of $V\subset\Ind_{\fg_0}^{\fg} V$ by $v$ 
and a generator of the trivial submodule of $\Ind_{\fg_0}^{\fg} V$
by $v'$. Retain notation of~\ref{descrAg}. 
Choose an ordered basis $\{ x_i\}_{i\in I}$
of $\fg_1$ in such a way that $[x_1,x_1]=0$. Write 
$v'=\sum_{J\subset I}c_Jx_J\otimes v$ where $c_J\in {\Bbb C}$.
One has $x_1x_J=0$ when
$1\in J$ and $x_1x_J=x_{\{1\}\cup J}$ otherwise. Since $x_1v'=0$, 
this implies that $c_{\emptyset}=0$.

Recall that $\Hom_{\fg} (\Ind_{\fg_0}^{\fg} V,\wdV)$ 
is one-dimensional and is spanned by an element $f$ such that
$f(V)\not=0$ and $f(\Ug\fg V)=0$. Thus $f(v')=0$ and so
the trivial submodule $\wdV$ is not a direct summand
of $\Ind_{\fg_0}^{\fg} V$.
\end{pf}

\subsubsection{}
Let $\fg$ be a basic classical Lie superalgebra which is 
neither simple Lie algebra nor $\osp(1,2l)$. Then $\fg_1$
contains a non-zero element $x$ such that $[x,x]=0$. 
From~(\ref{theta}) it follows that the multiplication by $T$
gives a $\fg$-map from the $\ad'\fg$-module generated 
by $1$ onto the $\ad\fg$-module generated by $T$. Since $T$
is a non-zero divisor, this map is an isomorphism. In particular, 
$T^2\in (\ad\Ug)T$ and $(\ad\Ug)T\cong\Ind_{\fg_0}^{\fg}V$
where $V$ is a trivial $\fg_0$-module. Using~\Lem{lemsepthm},
we conclude that $T^2$ spans a trivial $\ad\fg$-submodule of $\Ug$
which is not a direct summand. On the other hand, $1$
spans a trivial $\ad\fg$-submodule of $\Ug$ which is  a direct summand,
since $\Ug={\Bbb C}\oplus\Ug\fg$ as $\ad\fg$-module.

Assume now that Separation theorem holds for $\Ug$. Then 
$\Zg\cong (\cH)^{\fg}\otimes\Zg$ as $\Zg$-module. Therefore $(\cH)^{\fg}$
is one-dimensional. Thus all trivial $\ad\fg$-submodules of $\Ug$
are either direct summands or are not direct summands of $\Ug$.
This gives a contradiction.
  
\section{Questions}

\subsection{}
The centralizer of $\Ug_{\ol 0}$ contains $\tZg$. Do they coincide
provided that $\Lambda^{top}\fg_1$ is a trivial $\fg_0$-module?
Note that the condition on $\Lambda^{top}\fg_1$ is essential---
see~\Exa{exaW}.

\subsection{}
Let $C$ be the set of the elements of $\Ug$ which act
by a superconstant on each simple module. Clearly, $C$ is a subalgebra
of $\Ug$. By~\ref{evencase},
$C$ contains $\tZg$ if $\dim\fg_1$ is even.
Assume that $\dim\fg_1$ is even and that the
intersection of all graded primitive ideals of $\Ug$ is zero.
Does this imply that $C=\tZg$?
\subsection{}
In the case when $\fg$ is a basic classical Lie superalgebra both answers
are positive--- see~\Cor{supercaction}.

%%%%%%%%%%%%%%%  biblio.tex

\end{document}